\documentclass[12pt]{article}
\usepackage{amssymb}
\usepackage{amsmath}
\usepackage{latexsym}
\usepackage{mathrsfs}

\numberwithin{equation}{section}
\newtheorem{thm}[equation]{Theorem}

\newtheorem{defn}[equation]{Definition}
\newtheorem{rem}[equation]{Remark}
\newtheorem{lem}[equation]{Lemma}
\newtheorem{corol}[equation]{Corollary}

\title{Steklov-type eigenvalues associated with best Sobolev trace constants: domain perturbation and overdetermined systems.}

\author{Pier Domenico Lamberti}

\date{\ }

\begin{document}

\newcommand{\rea}{\mathbb{R}}

\maketitle
\centerline{\it  Dedicated to the  70th birthday of Professor Victor I. Burenkov}
\vspace{1cm}

%
%
%

\noindent
{\bf Abstract:}
We consider a variant of the classic Steklov eigenvalue problem, which arises  in the study of the best trace constant for functions
in  Sobolev space. We prove that the elementary symmetric functions of the eigenvalues depend real-analytically upon variation of the underlying
domain and we compute the corresponding Hadamard-type formulas for the shape derivatives. We also consider isovolumetric and isoperimetric
domain perturbations and we characterize the corresponding critical domains in terms of appropriate overdetermined systems. Finally, we prove that
balls are critical domains for the elementary symmetric functions of the eigenvalues subject to volume or perimeter constraint.
\\

\vspace{11pt}

\noindent
{\bf Keywords:}  Steklov  eigenvalues;
Laplace operator;  domain perturbation; overdetermined systems.

\vspace{6pt}
\noindent
{\bf 2000 Mathematics Subject Classification:} 35P15; 46E35; 35N25; 49Q12.

\section{Introduction}

Let $N\geq 2$ and $\Omega$ be a bounded domain ({\it i.e.}, a bounded connected open set) in ${\mathbb{R}}^N$ of class $C^2$. We consider the eigenvalue problem

\begin{equation}
\label{classic}
\left\{
\begin{array}{ll}
\Delta u=u, &\ {\rm in}\  \Omega ,\vspace{2mm}\\
\frac{\partial u}{\partial \nu }=\lambda u, &\ {\rm on }\  \partial \Omega ,
\end{array}\right.
\end{equation}
and we represent its eigenvalues by means of a divergent sequence
$$
0< \lambda_1[\Omega ]< \lambda_2[\Omega ] \le\dots\le\lambda_j[\Omega ]\le \dots
$$
where each eigenvalue is repeated according to multiplicity. Here $\nu $ denotes the outer unit normal to $\partial\Omega$. In this paper we study the dependence of $\lambda_j[\Omega ]$ on $\Omega $.

After the pioneering paper by Steklov~\cite{stek}, the homogeneous boundary condition $\partial u/\partial \nu =\lambda u$ is called Steklov boundary condition and it appears in the study of many boundary value problems.
In the literature, the Steklov condition is often imposed to harmonic functions in $\Omega $ and the eigenvalues of the problem
\begin{equation}
\label{steklovclassic}
\left\{
\begin{array}{ll}
\Delta u=0, &\  {\rm in}\ \Omega ,\vspace{2mm}\\
\frac{\partial u}{\partial \nu }=\lambda u, &\  {\rm on }\ \partial \Omega ,
\end{array}\right.
\end{equation}
are usually called Steklov eigenvalues. Problem (\ref{steklovclassic}) has important physical meaning and many properties of the Steklov eigenvalues are known. In particular, for $N=2$ problem (\ref{steklovclassic}) describes the vibration of a free membrane the total mass of which is uniformly distributed along the boundary. On the other hand, problem (\ref{classic}) seems to find its natural motivation in the frame of the Theory of Sobolev Spaces
and in particular in connection with the problem of traces. Recall that functions of the Sobolev space $W^{1,2}(\Omega )$ admit a trace in $L^2(\partial \Omega )$ and there exists a constant $C$ such that
\begin{equation}\label{trace}
\int _{\partial \Omega}u^2d\sigma \le C \int_{\Omega }u^2+|\nabla u|^2dx,
\end{equation}
for all $u\in W^{1,2}(\Omega )$.  Moreover, the trace operator from $W^{1,2}(\Omega )$ to $L^2(\partial \Omega )$ is compact.
By simple considerations, it turns out that
\begin{equation}\label{min}
\lambda_1[\Omega ]=\min_{\substack{u\in W^{1,2}(\Omega )\\ {\rm Tr}\, u\ne 0}} \frac{\int_{\Omega }|\nabla u|^2+|u|^2dx}{\int_{\partial \Omega}|u|^2d\sigma } >0,
\end{equation}
and that the minimizers are exactly the eigenfunctions corresponding to $\lambda_1[\Omega]$. In fact, problem (\ref{classic}) is a formulation of the Euler-Lagrange equation associated with (\ref{min}). Note that $\lambda_1[\Omega ]$ is the reciprocal of the best Sobolev trace constant $C$ for which (\ref{trace})
holds.

In the last decade problem (\ref{classic}) and its nonlinear generalizations have attracted the attention of several authors, see, {\it e.g.}, the rather widely quoted paper by Mart\'{\i}nez and Rossi~\cite{mar}. In particular,  Rossi~\cite{rossi} has conjectured that the the ball maximizes $\lambda_1[\Omega ]$ among all domains with fixed volume (this is known to be true for the first nontrivial eigenvalue of (\ref{steklovclassic}), see Weinstock~\cite{wei} and Brock~\cite{brock}, see also Henrot~\cite{hen}). In order to support his conjecture, Rossi~\cite{rossi}  has computed the shape derivative of $\lambda _1[\Omega ]$  and has proved that the ball is a critical point for $\lambda_1[\Omega ]$ under volume constraint.

In this paper, we develop the results of \cite{rossi} and we consider all eigenvalues $\lambda_j[\Omega]$. The main difficulty in dealing with higher eigenvalues is related to their multiplicity which may change when the domain is perturbed. This leads to complicated bifurcation phenomena, which clearly do not occur in the case of the simple eigenvalue $\lambda_1[\Omega ]$. The investigations carried out in the papers \cite{lala2004, lalacri, lamass}, have pointed out that in the case of multiple eigenvalues it is natural to consider the elementary symmetric functions of the  eigenvalues, which typically depend  real-analytically on the parameters involved in the problem. In this paper, we adopt the approach of \cite{lalacri}. Namely, we consider families of domains $\phi (\Omega )$ parameterized by  suitable diffeomorphisms $\phi $ of class $C^2$ defined on the fixed reference domain $\Omega $ and we prove that the elementary symmetric functions of the eigenvalues of problem (\ref{classic}) on $\phi (\Omega )$ depend real-analytically on $\phi $, see Theorem~\ref{sym}. Moreover, in formula (\ref{sym2}) we compute the Frech\'{e}t derivatives of such functions and we characterize the critical transformations $\phi $  subject to the volume constraint ${\rm Vol} (\Omega )= {\rm const}$ or the perimeter constraint ${\rm Per} (\phi (\Omega ))={\rm const}$, see Theorem~\ref{ov}. This leads to the formulation of the appropriate overdetermined systems (\ref{pr1}), (\ref{pr2}). It turns out that if $\phi (\Omega )$ is a ball, then $\phi $ is a critical point for the elementary symmetric functions of the eigenvalues with volume or perimeter constraint, and the corresponding overdetermined systems are satisfied.

We note that in the terminology of domain perturbation theory ({\it cf.} Henry~\cite{henry}), in this paper we adopt the `Lagrangian' point of view,
whilst the calculations in \cite{rossi} are performed in the `Eulerian' form.

\section{Notation and preliminaries}

In this paper the elements of ${\mathbb{R}}^N$ are thought as row vectors. The inverse of a matrix $A$ is denoted by $A^{-1}$ as opposed to the the inverse of a function $f$ which is denoted by  $f^{(-1)}$. The transpose of  a matrix $A$ is denoted by  $A^T$, and   the inverse of the transpose of $A$ is denoted by  $A^{-T}$.  The product of two  matrices $A$ and $B$ is denoted by $AB$. According to this notation,  $w_1w_2^T$ is the scalar product of two vectors  $w_1, w_2\in {\mathbb{R}}^N$.  The symbol $\nabla  $ denotes either the gradient of a real-valued function or
the Jacobian matrix of a vector-valued function. The differential of a function $F$ at a point $a$ is denoted by ${\rm d}F(a)$ and the value of ${\rm d}F(a)$ computed in $b$ is denoted by ${\rm d}F(a)[b]$.\\

  Let $\Omega$ be a bounded domain  in ${\mathbb{R}^N}$ of class $C^2$. Here and in the sequel it is understood that $N\geq 2$.
We consider the weak formulation of problem (\ref{classic})
\begin{equation}
\label{mainweak}
\int_{\Omega}\nabla u\nabla\varphi^T +  u\varphi \, dx = \lambda \int_{\partial \Omega}u\varphi \, d\sigma  ,\ \ {\rm for\ all }\ \varphi\in W^{1,2}(\Omega )
\end{equation}
in the unknowns $u\in W^{1,2}(\Omega )$ (the eigenfunction), $\lambda \in \mathbb{R}$ (the eigenvalue). Here $W^{1,2}(\Omega )$ denotes the Sobolev space of those functions
in $L^2(\Omega )$ with first order weak derivatives  in $L^2(\Omega )$ endowed with its standard norm (functions in $L^2(\Omega )$ are assumed to be real-valued) and $d\sigma$ denotes the $(N-1)$-dimensional surface measure in $\partial \Omega $.  Note that the integral in the right-hand side of (\ref{mainweak}) is well-defined since functions in the Sobolev space admit a trace in $L^2(\partial \Omega )$ in the sense that  there exists a  continuous linear operator ${\rm Tr}$ of $W^{1,2}(\Omega )$ to $L^2(\partial\Omega  )$ such that ${\rm Tr }\, u=u_{|\partial \Omega }$ for all $u\in W^{1,2}(\Omega )\cap C(\bar \Omega )$. Occasionally in the sequel we shall also use other standard Sobolev spaces $W^{m,p}(\Omega )$
of functions with weak derivatives up to order $m$ in $L^p(\Omega )$; we refer to Burenkov~\cite{bur} for the basic properties of these spaces.

We  now recall a standard procedure that enables us to reduce the study of problem (\ref{mainweak}) to the study of an eigenvalue problem
for a compact self-adjoint operator in the Hilbert space $W^{1,2}( \Omega )$ (which is  equipped with its standard scalar product defined by the left-hand side of (\ref{mainweak})).

We consider the Laplace operator $\Delta $ as an operator from $W^{1,2}(\Omega )$ to its dual $(W^{1,2}(\Omega ))'$
defined by $\Delta [u] [\varphi ]=-\int_{\Omega }\nabla u\nabla \varphi^T dx$ for all $u,\varphi \in W^{1,2}(\Omega )$ and we recall that the operator
$I-\Delta $ is an isometry of $W^{1,2}(\Omega )$ onto $(W^{1,2}(\Omega ))'$. Here $I$ is the natural embedding defined by $I[u][\varphi ]=\int_{\Omega }u\varphi dx$ for all $u,\varphi \in W^{1,2}(\Omega )$. Moreover, we consider the embedding $J$ of $L^2(\partial \Omega )$ into
$(W^{1,2}(\Omega ))'$ defined by $J[u][\varphi ]=\int_{\partial \Omega }u \varphi d\sigma $ for all $u \in L^2(\partial \Omega )$, $\varphi \in W^{1,2}(\Omega )$. Clearly, equation (\ref{mainweak}) is equivalent to equation
\begin{equation}
 (I-\Delta )^{-1}\circ J\circ {\rm Tr }\, [u]=\lambda^{-1} u
\end{equation}
in the unknown $u\in W^{1,2}( \Omega )$. Note that by using $\varphi =u$ as a test function in $(\ref{mainweak})$, it is immediate to see that any eigenvalue $\lambda $ is positive. Thus, we have the following

\begin{lem}\label{compact}
The operator $S_{\Omega }\equiv (I-\Delta )^{-1}\circ J\circ {\rm Tr }$ is a  compact non-negative self-adjoint operator in the Hilbert space $W^{1,2}(\Omega  )$. Moreover, $S_{\Omega }u=\mu u$ for $u\in W^{1,2}(\Omega )$, $\mu \in {\mathbb{R}}\setminus \{0\}$ if and only if $u$ is an eigenfunction of problem (\ref{mainweak}) corresponding to the eigenvalue $\lambda =\mu^{-1}$.
\end{lem}

{\bf Proof.}  By equality
$$
<S_{\Omega }u_1,u_2>_{W^{1,2}( \Omega )}= (I-\Delta)[(I-\Delta)^{-1} \circ J\circ {\rm Tr }\, [u_1]  ][u_2]=J\circ {\rm Tr}\, [u_1 ][u_2],
$$
for all $u_1,u_2\in W^{1,2}(\Omega  )$ and by the symmetry of the pairing $J\circ {\rm Tr}\, [\cdot  ][\cdot ]$ we immediately deduce that $S_{\Omega }$ is a non-negative self-adjoint operator. The compactness of $S_{\Omega }$ follows by the compactness of the trace operator ${\rm Tr}$ (see,  {\it e.g.}, Kufner, John, Fu\v{c}\'{\i}k~\cite[Thm.~6.10.5]{kujo}).
The proof of the last part of the statement  is straightforward.\hfill $\Box$\\

It is easy to see that ${\rm Ker } S_{\Omega }= W^{1,2}_0(\Omega )$, where $W^{1,2}_0(\Omega )$ is the closure in $W^{1,2}(\Omega )$ of the space of functions of class $C^{\infty }$ with compact support in $\Omega $.
Since the operator $S_{\Omega }$ is compact and self-adjoint in the Hilbert space $W^{1,2}(\Omega)$ and its kernel has infinite codimension in $W^{1,2}(\Omega )$, the spectrum  of  $\sigma (S_{\Omega })$ is discrete and $\sigma (S_{\Omega })\setminus \{0\}$  consists of eigenvalues of finite multiplicity which can be represented as a non-increasing sequence
$\mu_j[S_{\Omega }]$, $j\in {\mathbb{N}}$ converging to zero (as customary, we repeat each eigenvalue as many times as its multiplicity). Thus we have the following

\begin{corol} Let $\Omega$ be a bounded domain in ${\mathbb{R}}^N$ of class $C^2$. The eigenvalues of problem (\ref{mainweak}) are positive, have finite multiplicity and can be represented as a  non-decreasing divergent sequence $\lambda_j[\Omega ]$, $j\in {\mathbb{N}}$ where each eigenvalue is repeated according to its multiplicity. Moreover,
\begin{equation}
\label{minmax}
\lambda_j[\Omega ]=\min_{\substack{E\subset W^{1,2}(\Omega ) \\  E\cap  W^{1,2}_0(\Omega )=\{0\}\\    {\rm dim } E=j}} \max_{\substack{ u\in E \\  u\ne 0}} \frac{\int_{\Omega }|\nabla u|^2+|u|^2dx}{\int_{\partial \Omega}|u|^2d\sigma },
\end{equation}
for all $j\in {\mathbb{N}}$ and formula (\ref{min}) holds. Furthermore, $\lambda_1[\Omega ]$ is simple and its eigenfunctions do not change sign in $\Omega$.
\end{corol}

{\bf Proof.} The proof of the first part of the statement can be easily deduced by Lemma~\ref{compact} and standard spectral theory and by recalling that $\lambda_j[\Omega ]=\mu_j^{-1}[S_{\Omega } ]$ for all $j\in {\mathbb{N}}$.  For the second part we proceed exactly as in the well-known case of the Dirichlet Laplacian. Namely, let $u$ be a nonzero eigenfunction associated with $\lambda_1[\Omega ]$. Since $u$ is a minimizer in (\ref{min}) then also $|u|$ is a minimizer in (\ref{min}). Thus, $| u|$ is an eigenfunction associate with $\lambda_1[\Omega]$ and by the Harnack inequality $| u|$ cannot vanish in $\Omega$. Hence $u$ does not change sign in $\Omega$.  Finally, given two nonzero eigenfunctions  $u_1, u_2$ associated with  $\lambda_1[\Omega ]$ and $c\in {\mathbb{R}}$ such that
$\int_{\Omega }u_1-cu_2dx=0$, it follows that the eigenfunction $u_1-cu_2$ must be identically zero, hence $u_1=cu_2$.\hfill $\Box$\\

\begin{rem}
Another way of reducing problem (\ref{mainweak}) to an eigenvalue problem for a compact non-negative self-adjoint operator in Hilbert space is to
consider the operator ${\rm Tr}\circ (I-\Delta )^{-1}\circ J$ in  $L^2(\partial \Omega )$. This operator may be used instead of $S_{\Omega }$ in all our arguments below without any essential changes. With regard to this, we note that, in order  to normalize  eigenfunctions, the scalar product of $L^2(\partial \Omega )$  would  be more natural than the scalar product of $W^{1,2}(\Omega )$ (see also Remark~\ref{remnorm} in Section~\ref{domsec}).  However, the relation between the eigenfunctions of ${\rm Tr}\circ (I-\Delta )^{-1}\circ J$ and those of problem (\ref{mainweak}) is  slightly involved. Thus, we prefer to use the operator $S_{\Omega }$ since its eigenfunctions are functions defined on $\Omega $ and coincide with the eigenfunctions of problem (\ref{mainweak}).
\end{rem}

\section{Domain perturbation}

\label{domsec}

Let $\Omega$ be a fixed bounded domain  in ${\mathbb{R}}^N$ of class $C^2$. We consider problem (\ref{mainweak}) on a class of domains diffeomorphic
to $\Omega $. Namely, we consider the class of domains $\{\phi (\Omega ):\ \phi \in {\mathcal{A} }_{\Omega }\} $ where ${\mathcal {A}}_{\Omega }$ is defined by
\begin{equation}
{\mathcal{A}}_{\Omega }=\left\{\phi \in C^2(\bar \Omega\, ; {\mathbb{R}}^N):\ \phi \ {\rm is \ injective\ and}\  {\rm det }\nabla \phi (x)\ne 0, \ {\rm for\ all}\ x\in \bar\Omega    \right\},
\end{equation}
and $ C^2(\bar \Omega\, ; {\mathbb{R}}^N)$ is the space of the functions from $\bar \Omega $ to ${\mathbb{R}}^N$ of class $C^2$.
Note that if $\phi \in {\mathcal {A}}_{\Omega }$ then $\phi (\Omega )$ is a bounded domain of class $C^2$ in ${\mathbb{R}}^N$ and that $\partial \phi(\Omega )=\phi (\partial \Omega)$. Thus, it is possible to set the eigenvalue problem (\ref{mainweak}) on $\phi (\Omega )$ for all $\phi \in {\mathcal{A}}_{\Omega }$ and to consider the eigenvalues $\lambda_j[\phi (\Omega )]$ as functions of $\phi $. For simplicity, we set
$\lambda_j[\phi ]=\lambda_j[\phi (\Omega )]$.

In this section, we study the dependence of $\lambda_j[\phi ]$ upon variation of $\phi \in {\mathcal{A}}_{\Omega }$. It is understood that the space $C^2(\bar \Omega\, ; {\mathbb{R}}^N)$ is endowed with the $C^2$-norm defined by
$ \| \phi \| =\sum_{i=1}^N \sum_{|\alpha |\le 2}\| D^{\alpha }\phi_i  \|_{\infty } $ for all $\phi =(\phi_1, \dots ,\phi_N)\in C^2(\bar \Omega\, ; {\mathbb{R}}^N) $, where $\|\cdot \|_{\infty }$ is the standard sup-norm. It turns out that  the set ${\mathcal{A}}_{\Omega }$ is an open set in the space $C^2(\bar \Omega\, ; {\mathbb{R}}^N)$ (see \cite[Lemma~2.2]{lalagla}). Thus, it makes sense to study differentiability and analyticity properties of the maps $\lambda_j[\cdot ]$ defined on ${\mathcal{A}}_{\Omega }$.

Following \cite{lala2004}, we fix a finite set of indexes $F\subset \mathbb{N}$
and we consider those maps $\phi\in {\mathcal{A}}_{\Omega }$ for which the eigenvalues
with indexes in $F$  do not coincide with eigenvalues with indexes not
in $F$; namely we set
$$
{\mathcal { A}}_{\Omega }[F]= \left\{\phi \in {\mathcal { A}}_{\Omega }:\
\lambda_j[\phi ]\ne \lambda_l[\phi],\ \forall\  j\in F,\,   l\in \mathbb{N}\setminus F
\right\}
$$
and
$$
\Theta_{\Omega } [F]= \left\{\phi\in {\mathcal { A}}_{\Omega }[F]:\ \lambda_{j_1}[\phi ]
=\lambda_{j_2}[\phi ],\, \
\forall\ j_1,j_2\in F  \right\} .
$$

Then we can prove the following theorem, where $\nabla_Tv (y)$ denotes the orthogonal projection of $\nabla v (y) $ to the tangent plane
of $\partial\phi (\Omega ) $ at the point $y\in \partial \phi (\Omega )$. Note that by standard elliptic  theory ({\it cf.} Agmon, Douglis and Nirenberg~\cite{agmon}; see also the recent paper by Castro Triana~\cite{castro}), any eigenfunction $v$ in $\phi (\Omega )$ is of class $W^{2,2}(\phi (\Omega ))$ hence $\nabla v$ admits a trace in $L^2(\partial \phi (\Omega ))$. Thus, $\nabla_Tv$ and  the normal derivative  $\partial v /\partial\nu  $  are well-defined and $|\nabla v|^2=|\nabla_Tv|^2+|\partial v /\partial\nu  |^2$. Although it is not necessary, it may be useful to recall that  $C^{1,\alpha }$-type results ({\it cf.},  {\it e.g.}, Lieberman~\cite{liebe}) imply that the gradient of $v$ is in fact continuous up to the boundary of $\phi (\Omega )$, hence $\nabla_T v$ and $\partial v/ \partial\nu   $ are actually defined everywhere in $\partial \phi (\Omega )$.

\begin{thm}
\label{sym}
Let $\Omega$ be a bounded domain in $\mathbb{R}^N$ of class $C^2$ and $F$ a finite subset of ${\mathbb{N}}$.
Then  ${\mathcal{A}}_{\Omega }[F]$ is an open set in $  C^2(\bar \Omega\, ; {\mathbb{R}}^N)  $ and the
elementary symmetric functions
\begin{equation}
\label{sym1}
\Lambda_{F,h}[\phi ]=\sum_{ \substack{ j_1,\dots ,j_h\in F\\ j_1<\dots <j_h} }
\lambda_{j_1}[\phi ]\cdots \lambda_{j_h}[\phi ],\ \ \ h=1,\dots , |F|
\end{equation}
are real-analytic in ${\mathcal { A}}_{\Omega }[F]$. Moreover,
if $\phi\in \Theta_{\Omega } [F]$ is such that the eigenvalues $\lambda_j[\phi ]$ assume the common value $\lambda_F[\phi ]$
for all $j\in F$, then the differential of the functions $\Lambda_{F,h}$
at the point $\phi$
is given
by the formula
\begin{eqnarray}
\label{sym2}\lefteqn{
{\rm d}\Lambda_{F,h}(\phi)[\psi] =
\lambda_{F}^{h}[\phi]
{
|F|-1
\choose
h-1
}
\sum_{l\in F}
\int_{\partial \phi (\Omega )} \left[  |\nabla_T v_l|^2 \nonumber \right. }\\ & &\left. \qquad\qquad\qquad\quad\quad\quad+(1-\lambda_F[\phi]H-\lambda_F^2[\phi ])v_l^2\right] (\psi \circ \phi ^{(-1)})\nu^T
 d\sigma ,
\end{eqnarray}
for all $\psi \in C^2(\bar \Omega\, ; {\mathbb{R}}^N) $,
where $\{ v_l\}_{l\in F}$ is an orthonormal basis in $W^{1,2} ( \phi (\Omega ))$ of the eigenspace corresponding to $\lambda_F[\phi ]$,
 and $H={\rm div }\nu $ is the mean curvature of  $\partial \phi (\Omega)$ (the sum of the principal curvatures).
\end{thm}

\begin{rem}\label{remnorm} By looking at the weak formulation (\ref{mainweak}), it is immediate to realize that if one considers  eigenfunctions $v_l$  normalized with respect to the scalar product of $L^2(\partial \phi (\Omega ))$ then the factor $\lambda_F^{h}[\phi ]$ in formula (\ref{sym2}) should be replaced by $\lambda_F^{h-1}[\phi ]$. In particular it follows that formula (\ref{sym2}) computed for  $|F|=h=1$ agrees (up to a misprint) with the formula in Rossi~\cite[(1.5)]{rossi} concerning the first eigenvalue.
\end{rem}

The rest of this Section is devoted to proving Theorem~\ref{sym}. We note that the proof of Theorem~\ref{sym} is based on a general  result in \cite{lala2004} concerning families of  compact self-adjoint operators in Hilbert space with variable scalar product. In order to apply that result we need to consider the eigenvalue problem (\ref{mainweak}) on the variable domain $\phi (\Omega )$ and pull it back
to the fixed domain $\Omega$ by means of a change of variables. In particular, for a fixed $\phi \in {\mathcal{A}}_{\Omega }$, it is convenient to consider the operators $\Delta$, $I$, and $J$ on $\phi (\Omega )$ defined in Section 2, and  pull them back to $\Omega $. In this way we obtain the operators
$\Delta_{\phi }$, $I_{\phi }$ of $W^{1,2}(\Omega )$ to $(W^{1,2}(\Omega ))'$ and $J_{\phi }$ of $L^2(\partial \Omega )$  to $(W^{1,2}(\Omega ))'$ defined as follows:
\begin{equation}
\label{laplacephi1}
\Delta_{\phi}\left[u\right]\left[\varphi\right]=\Delta\left[u\circ\phi^{\left(-1\right)}\right]\left[\varphi \circ\phi^{\left(-1\right)}\right],\ I_{\phi }[u][\varphi ]=I\left[u\circ\phi^{\left(-1\right)}\right]\left[\varphi \circ\phi^{\left(-1\right)}\right],
\end{equation}
for all $u,\varphi \in W^{1,2}(\Omega )$, and
\begin{equation}
J_{\phi }[u][\varphi ]=J\left[u\circ\phi^{\left(-1\right)}\right]\left[\varphi \circ\phi^{\left(-1\right)}\right] ,
\end{equation}
for all $u\in L^2(\partial \Omega )$, $\varphi\in W^{1,2}(\Omega )$.
By means of standard  calculus
and a change of variables one can check that
\begin{equation}
\label{laplacephi2}
\Delta_{\phi}\left[u\right]\left[\varphi\right]=
-\int_{\Omega}\nabla u \left(\nabla \phi\right)^{-1}\left(\nabla \phi\right)^{-T}\ \nabla \varphi ^T\left|\mathrm{det}\nabla \phi\right|dx,
\end{equation}
\begin{equation}\label{laplacephi2bis}
I_{\phi }[u ][\varphi] = \int_{\Omega }u\varphi \left|\mathrm{det}\nabla \phi\right|dx ,
\end{equation}
for all $u,\varphi \in W^{1,2}(\Omega )$,
and
\begin{equation}\label{bdchange}
J_{\phi }[u][\varphi ]=\int_{\partial \Omega }u\varphi \left|\nu (\nabla \phi )^{-1}   \right| \left|\mathrm{det}\nabla \phi\right|d\sigma ,
\end{equation}
$u\in L^2(\partial \Omega )$, $\varphi \in W^{1,2}(\Omega )$.

Formulas (\ref{laplacephi2}) and (\ref{laplacephi2bis}) suggest introducing in $W^{1,2}( \Omega )$ the following scalar product
\begin{equation}\label{scalarphi}
<u_1,u_2 >_{\phi }=  \int_{\Omega}(  \nabla u_1 \left(\nabla \phi\right)^{-1}\left(\nabla \phi\right)^{-T}\ \nabla u_2 ^T +u_1u_2  ) \left|\mathrm{det}\nabla \phi\right|dx,
\end{equation}
for all $u_1,u_2 \in W^{1,2}(\Omega )$.  We denote by $W^{1,2}_{\phi}( \Omega )$ the space $W^{1,2}(\Omega  )$ equipped with the
scalar product defined by (\ref{scalarphi}). Note that the scalar product (\ref{scalarphi}) is equivalent to the standard scalar product of
$W^{1,2}(\Omega  )$, hence $  W^{1,2}_{\phi}( \Omega )$ is in fact a Hilbert space. Note also that the map $C_{\phi}$ of $W^{1,2}(\phi (\Omega ))$ to
$W^{1,2}(\Omega ) $ defined by $C_{\phi }[v]=v\circ \phi $ for all $v\in W^{1,2}_{\phi }(\phi(\Omega ))$ is an isometry.
By setting
$$
T_{\phi }= (I_{\phi }-\Delta_{\phi })^{-1}\circ J_{\phi}\circ {\rm Tr} ,
$$
it is immediate to verify the validity of the following

\begin{lem}\label{isolem}
The operator $T_{\phi }$ of $W^{1,2}_{\phi }( \Omega )$ to itself satisfies
$$
T_{\phi }= C_{\phi }\circ  S_{\phi (\Omega )} \circ C_{\phi ^{-1}} ,
$$
hence it is unitarily equivalent to $S_{\phi (\Omega )}$. In particular, $T_{\phi }$ is a compact non-negative self-adjoint operator
in  $W^{1,2}_{\phi }( \Omega )$ the eigenvalues of which coincide with the eigenvalues of $S_{\phi (\Omega )} $. Moreover,  $v\in W^{1,2}(\phi (\Omega ))$ is an eigenfunction of $S_{\phi (\Omega )} $ if and only if $v\circ \phi$ is an eigenfunction of $T_{\phi }$.
\end{lem}

In order to prove formula (\ref{sym2}) we need the following technical lemma. To avoid heavy notation, we use the same symbol
$\nu $ to denote the  outer unit normal to the boundaries of different domains (the context is usually clear enough to understand which domain $\nu $ refers to; for example, $\nu $ refers to $\partial \Omega $ in (\ref{techlem0}) and  $\partial \phi (\Omega )$ in (\ref{techlem1})).

\begin{lem}\label{techlem}
Let $\Omega $ be a bounded domain of class $C^2$. Let $u\in W^{2,1}(\Omega)$ be fixed. The function $B$ of ${\mathcal {A}}_{\Omega }$ to ${\mathbb{R}}$ defined by
\begin{equation}
\label{techlem0} B[\phi ] =\int_{\partial \Omega }u \left|\nu (\nabla \phi )^{-1}   \right| \left|\mathrm{det}\nabla \phi\right|d\sigma ,\end{equation}
for all  $\phi\in {\mathcal{A}}_{\Omega }$, is real-analytic and the differential at the point $\phi \in {\mathcal{A}}_{\Omega }$ is given by
\begin{equation}\label{techlem1}
{\rm d} B(\phi )[\psi]=\int_{\partial \phi (\Omega )} \biggl( Hv+\frac{\partial{v}}{\partial \nu }\biggr)(\psi\circ \phi^{(-1)})\nu^T d\sigma -\int_{\partial\phi( \Omega )}(\psi\circ \phi^{(-1)})\nabla v^Td\sigma ,
\end{equation}
for all $\psi\in C^2(\Omega\, ; {\mathbb{R}}^N)$, where $v=u\circ \phi^{(-1)}$.
\end{lem}

{\bf Proof}
Since ${\rm det} \nabla\phi $ does not vanish on $\partial\Omega $, the weight  $\left|\nu (\nabla \phi )^{-1}   \right| | {\rm det} \nabla \phi|$ does not vanish on $\partial \Omega $. Thus, by a continuity argument $\inf_{\partial \Omega } \left|\nu(\nabla \phi )^{-1}   \right| | {\rm det} \nabla \phi| $ is strictly positive in a suitable neighborhood in ${\mathcal{A}}_{\Omega }$ of any fixed $\tilde \phi \in {\mathcal{A}}_{\Omega }$. Thus, the map $B$ is real-analytic since it is  the composition of real-analytic maps.

 By standard
calculus, it is easy to see that differentiating ${\mathrm{det}}\nabla\phi$ with respect to $\phi $ yields \begin{equation}
\label{der4}
\left[
\left({\rm d} \left({\mathrm{det}}\nabla\phi\right)[\psi]\right)\circ
\phi^{(-1)}\right]{\mathrm{det}}\nabla\phi^{(-1)}=
{\mathrm{div}}\left(\psi\circ \phi^{(-1)} \right)\, ,
\end{equation}
for all $\psi\in C^2(\Omega\, ; {\mathbb{R}}^N)$.
Moreover
\begin{equation}
\label{bdr1}
{\rm d} |\nu (\nabla\phi) ^{-1} |= -\frac{\nu (\nabla\phi)^{-1}[ \nabla\psi (\nabla\phi )^{-1} +(\nabla\phi )^{-T}(\nabla\psi)^T](\nabla\phi)^{-T}\nu^T  }{2|\nu (\nabla\phi )^{-1}|}.
\end{equation}
We note that the outer unit normal to $\partial \phi (\Omega )$ at the point $\phi (x)\in \partial \phi (\Omega )$ is given up to the sign by
\begin{equation}
\label{bdr2}
\frac{\nu (x)  (\nabla\phi (x))^{-1}}{|\nu (x)(\nabla\phi (x))^{-1}|} ,
\end{equation}
where $\nu (x)$ is the  outer unit normal to $\partial \Omega $ at the point $x\in \partial \Omega $.
Thus by combining (\ref{bdr1}), (\ref{bdr2}) and by means of a change of variable (see also formula (\ref{bdchange})) and by the Divergence Theorem applied in $\phi (\Omega )$ we have that
\begin{eqnarray}\label{bdr3}\lefteqn{\int_{\partial\Omega }u {\rm d} |\nu (\nabla\phi)^{-1} |[\psi ] |{\rm det }\nabla\phi |d\sigma } \nonumber \\ & & \qquad\qquad =
-\frac{1}{2}\int_{\partial \phi (\Omega )}v\nu [\nabla(\psi\circ\phi^{(-1)})+\nabla(\psi\circ\phi^{(-1)})^T]\nu^Td\sigma \nonumber \\
& & \qquad\qquad
=-\int_{\phi (\Omega )} {\rm div }(v \nabla (\psi\circ\phi^{(-1)} )\nu^T )dy\, .
\end{eqnarray}
In (\ref{bdr3}) and in the sequel it is understood that the unit vector field $\nu (y)$
defined for all $y\in \partial \phi (\Omega ) $ is extended as a unit vector field  of class $C^1$ in a neighborhood of $\partial \phi (\Omega )$ and then extended in any manner as a vector field of class $C^1$  defined in the whole of $\phi (\Omega )$ (see also Henry~\cite[p.~15]{henry}).
By (\ref{der4}) and (\ref{bdr3}) it follows that
\begin{eqnarray}\label{bdr4}
{\rm d}B[\phi][\psi]= \int_{\partial \phi (\Omega )}v{\rm div}(\psi \circ \phi^{(-1)})d\sigma -\int_{\phi (\Omega )}{\rm div }(v \nabla (\psi \circ \phi^{(-1)})   \nu^T )dy\nonumber \\
= \int_{\partial \phi (\Omega )}v{\rm div}(\psi \circ \phi^{(-1)})d\sigma -\int_{\phi (\Omega )}\nabla v \nabla(\psi \circ \phi ^{(-1)} )\nu ^Tdy\nonumber \\
-\int_{\phi (\Omega )}v{\rm div }(\nabla(\psi \circ \phi ^{(-1)})\nu^T )dy\, .\end{eqnarray}
To proceed with these computations we need  to differentiate the vector field $\nu $ twice. Thus, since $\nu $ is of class $C^1$, we approximate $\nu$ by a sequence of $C^2$ vector fields $\nu ^{(n)}=(\nu^{(n)}_1,\dots ,\nu^{(n)}_N )$, $n\in {\mathbb{N}}$ such that $\nu ^{(n)}\to \nu $ in $C^1({\overline{\phi (\Omega )}})$ as $n\to \infty$. To shorten our notation we set $\zeta =\psi \circ \phi ^{(-1)} $. By applying repeatedly the Divergence Theorem, we obtain
\begin{eqnarray}
\lefteqn{\int_{\phi (\Omega )} \nabla v\nabla\zeta (\nu^{(n)})^Tdy= \sum_{r,s=1}^N\int_{ \phi (\Omega )}\frac{\partial \zeta_r}{\partial y_s}\frac{\partial v}{\partial y_r}\nu^{(n)}_sdy= \int_{\partial \phi (\Omega )}(\nu^{(n)}\nu^T)  \nabla v \zeta ^Td\sigma }\nonumber\\
& &\quad -\sum_{r,s=1}^N  \int_{\phi (\Omega )} \zeta_r\frac{\partial^2v}{\partial y_r\partial y_s}\nu_s^{(n)}+\zeta_r\frac{\partial v}{\partial y_r}\frac{\partial \nu ^{(n)}_s}{\partial y_s}dy=\int_{\partial \phi (\Omega )} (\nu^{(n)}\nu^T) \nabla v\zeta^Td\sigma \nonumber \\
& &\quad -\sum_{r,s=1}^N\int_{\phi (\Omega )}   \zeta_r\frac{\partial^2v}{\partial y_r\partial y_s}\nu_s^{(n)}dy-\sum_{s=1}^N\int_{\partial \phi (\Omega )}\frac{\partial \nu^{(n)}_s}{\partial y_s} v \nu\zeta^Td\sigma  \nonumber \\
& &\quad \label{bdr5} +\sum_{r,s=1}^N\int_{\phi (\Omega )}v\frac{\partial \zeta_r}{\partial y_r}\frac{\partial \nu_s^{(n)}}{\partial y_s} + v\zeta_r\frac{\partial^2\nu_s^{(n)}}{\partial y_r\partial y_s}dy\, ,
\end{eqnarray}
and
\begin{eqnarray}\label{bdr6}
\lefteqn{\int_{\phi (\Omega )}v{\rm div }(\nabla\zeta (\nu^{(n)})^T)dy =\sum_{r,s=1}^N\int_{\phi (\Omega )}v \frac{\partial ^2\zeta_r}{\partial y_s\partial y_r}\nu_s^{(n)}+v\frac{\partial \zeta_r}{\partial y_s}\frac{\partial \nu_s^{(n)}}{\partial y_r}dy}\nonumber  \\
& & \qquad
=\int_{\partial \phi (\Omega )} (\nu^{(n)}\nu^T)   v{\rm div \zeta }d\sigma -\sum_{r,s=1}^N\int_{\phi (\Omega )}\frac{\partial v}{\partial y_s}\frac{\partial \zeta_r}{\partial y_r}\nu_s^{(n)}+v\frac{\partial \zeta_r}{\partial y_r}\frac{\partial \nu_s^{(n)}}{\partial y_s}dy\nonumber\\
& & \qquad+\sum_{r,s=1}^N\int_{\partial \phi (\Omega )}v\zeta _r\nu_s\frac{\partial \nu_s^{(n)}}{\partial y_r}d\sigma -\sum_{r,s=1}^N\int_{\phi (\Omega )}v\zeta _r\frac{\partial^2\nu_s^{(n)} }{\partial y_r\partial y_s}+\zeta_r\frac{\partial \nu_s^{(n)}}{\partial y_r}\frac{\partial v}{\partial y_s}dy\nonumber\\
& & \qquad=\int_{\partial \phi (\Omega )} (\nu^{(n)}\nu^T) v{\rm div \zeta }d\sigma -\sum_{r,s=1}^N\int_{\phi (\Omega )}\frac{\partial v}{\partial y_s}\frac{\partial \zeta_r}{\partial y_r}\nu_s^{(n)}+v\frac{\partial \zeta_r}{\partial y_r}\frac{\partial \nu_s^{(n)}}{\partial y_s}dy\nonumber \\
& & \qquad+\sum_{r,s=1}^N\int_{\partial \phi (\Omega )}v\zeta _r\nu_s\frac{\partial \nu_s^{(n)}}{\partial y_r}d\sigma -\sum_{r,s=1}^N\int_{\phi (\Omega )}v\zeta _r\frac{\partial^2\nu_s^{(n)} }{\partial y_r\partial y_s}- \zeta_r\frac{\partial^2 v}{\partial y_r\partial y_s}\nu_s^{(n)}\nonumber\\
& & \qquad-\frac{\partial \zeta_r}{\partial y_r}\nu_s^{(n)}\frac{\partial v}{\partial y_s}dy -\int_{\partial \phi(\Omega )} (\nu^{(n)} \nabla v^T )    \zeta \nu  ^Td\sigma\, .
\end{eqnarray}

Thus
\begin{eqnarray}\label{bdr7}& &
\int_{\phi (\Omega )} \nabla v\nabla\zeta (\nu^{(n)})^Tdy+\int_{\phi (\Omega )}v{\rm div }(\nabla\zeta (\nu^{(n)})^T)dy= \int_{\partial \phi (\Omega )} (\nu^{(n)}\nu ^T) v{\rm div \zeta }d\sigma\nonumber \\
& &\qquad\qquad  +\int_{\partial \phi (\Omega )}(\nu^{(n)}\nu ^T)  \nabla v\zeta^Td\sigma
-\sum_{s=1}^N\int_{\partial \phi (\Omega )}\frac{\partial \nu^{(n)}_s}{\partial y_s} v\nu \zeta^Td\sigma\nonumber  \\
& &\qquad\qquad  +\sum_{r,s=1}^N\int_{\partial \phi (\Omega )}v\zeta _r\nu_s\frac{\partial \nu_s^{(n)}}{\partial y_r}d\sigma
-\int_{\partial \phi(\Omega )}    ( \nu^{(n)}  \nabla v^T )    \zeta \nu  ^T             d\sigma\, .
\end{eqnarray}

By passing to the limit in (\ref{bdr7}) as $n\to \infty$,  by observing that $\sum_{s=1}^N\nu_s\frac{\partial \nu_s}{\partial y_r}=0$ on $\partial \phi (\Omega )$ for all $r=1, \dots , N$, and by recalling that $H={\rm div }\nu $ we obtain that

\begin{eqnarray}\label{bdr8}\lefteqn{
\int_{\phi (\Omega )} \nabla v\nabla\zeta \nu^Tdy+\int_{\phi (\Omega )}v{\rm div }(\nabla\zeta \nu^T)dy=\int_{\partial \phi (\Omega )}v{\rm div \zeta }d\sigma}\nonumber \\
& &\qquad  + \int_{\partial \phi (\Omega )}\nabla v\zeta^Td\sigma
-\int_{\partial \phi (\Omega )} H v \nu \zeta^Td\sigma -\int_{\partial \phi(\Omega )}\frac{\partial v}{\partial \nu }   \zeta\nu^Td\sigma\, .
\end{eqnarray}

By combining (\ref{bdr4}) and (\ref{bdr8}), we conclude.\hfill $\Box$\\

{\bf Proof of Theorem~\ref{sym}:}\quad We denote by $B_s((W^{1,2}(\Omega ))^2; \mathbb{R})$ the Banach space of the
bounded real-valued bilinear symmetric forms in $W^{1,2}(\Omega )$ endowed with its usual
norm. We denote by ${\mathcal {Q}}((W^{1,2}(\Omega ))^2; \mathbb{R} )  $ the
open subset of  $B_s((W^{1,2}(\Omega ))^2; \mathbb{R})$  of those bilinear
forms which define a scalar product in $ W^{1,2}(\Omega )$ topologically equivalent
to the standard scalar product of $W^{1,2}(\Omega ) $.
Finally we denote by $K(  W^{1,2}(\Omega ); W^{1,2}(\Omega ) )  $  the Banach space
of the compact linear operators in $ W^{1,2}(\Omega )$  endowed with its usual norm.
We note that the map from ${\mathcal{A}}_{\Omega } $ to
${\mathcal {Q}}((W^{1,2}(\Omega ))^2; {\mathbb{R}} )$ $ \times
K(W^{1,2}(\Omega );W^{1,2}(\Omega ) ) $ which takes $\phi\in {\mathcal{A}}_{\Omega } $
to the couple $(<\cdot ,\cdot >_{\phi}, T_{\phi })$ is real-analytic, since
it is the composition of real-analytic functions (see also the beginning of the proof of Lemma~\ref{techlem}). Moreover, by Lemma~\ref{isolem} $T_{\phi }$
is self-adjoint with respect to the scalar product
$<\cdot ,\cdot >_{\phi} $ and the eigenvalues $\mu_j[\phi] $ of
$T_{\phi }$ coincide with the reciprocals of $\lambda_j[\phi ]$.
Clearly, the set  ${\mathcal { A}}_{\Omega }[F]$  coincides
with the set
$\{\phi\in {\mathcal{A}}_{\Omega }:\
\mu_j[\phi]\ne \mu_l[\phi],\  \forall j\in F,\   l\in \mathbb{N}\setminus F
\} $.
Thus, by applying \cite[Thm.~2.30]{lala2004},
it follows that
${\mathcal { A}}_{\Omega }[F]$ is an open set in $C^2(\bar \Omega\, ; {\mathbb{R}}^N)$ and the
functions which take $\phi\in {\mathcal { A}}_{\Omega }[F] $  to
\begin{equation}
\label{sym3}
\Gamma_{F,h}[\phi ]=\sum_{\substack{j_1,\dots ,j_h\in F\\ j_1<\dots <j_h}}
\mu_{j_1}[\phi ]\cdots \mu_{j_h}[\phi ]
\end{equation}
are real-analytic for all $h=1, \dots , |F|$.
Since
\begin{equation}
\label{sym4}
\Lambda_{F,h}[\phi ] = \frac{\Gamma_{F,|F|-h}[\phi ] }{\Gamma_{F,|F|}[\phi ]},
\end{equation}
for all $h=1, \dots ,|F|$, where $\Gamma_{F,0}[\phi ]= 1  $,
it follows that $\Lambda_{F,h}[\phi ]$ depends real-analytically
on $\phi\in {\mathcal { A}}_{\Omega }[F]$.

It remains to  prove formula (\ref{sym2}).  Let $\phi \in \Theta_{\Omega }[F]$, $\lambda_F[\phi ]$ and $\{ v_l\}_{l\in F}$ be as in the statement. We set $u_l=v_l\circ \phi $ for all $l\in F$ and we note that $\{ u_l\}_{l\in F}  $ is an orthonormal basis in $W^{1,2}_{\phi }(\Omega )$ for the eigenspace corresponding to the eigenvalue $\lambda ^{-1}_F[\phi ]$ of the operator $T_{\phi }$.
By
\cite[Thm.~2.30]{lala2004}, it follows that
\begin{equation}
\label{sym5}
{\rm d}\Gamma_{F,h}(\phi )[\psi]=
\lambda_{F}^{1-h}[\phi]    {|F|-1\choose h-1}
\sum_{l\in F}
<{\rm d}T_{\phi}[\psi][u_l], u_l>_{\phi},
\end{equation}
for all $\phi\in {\mathcal { A}}_{\Omega }[F] $, $\psi\in C^2(\bar \Omega\, ;{\mathbb{R}}^N) $.
Thus we have to compute the summands in the right-hand side of (\ref{sym5}). For the sake of clarity, we believe it is
better to compute $<{\rm d}T_{\phi}[\psi][u_l], u_m>_{\phi}$ for all $l,m\in F$.

To shorten our notation we set $\lambda =\lambda_F[\phi ]$ and  $\zeta =\psi \circ\phi ^{(-1)} $.
Since $T_{\phi }u_l=\lambda^{-1}u_l  $ for all $l\in F$,  we obtain
\begin{eqnarray}\label{dedelta0}\lefteqn{
<{\rm d}T_{\phi}[\psi][u_l], u_m>_{\phi}}\nonumber  \\& &  =<(I_{\phi }-\Delta_{\phi }  )^{(-1)}\circ {\rm d}J_{\phi }[\psi ][{\rm Tr}\, u_l] + {\rm d}(  I_{\phi }-\Delta_{\phi }  )^{(-1)}[\psi ]\circ J_{\phi }[{\rm Tr}\, u_l], u_m>_{\phi }\nonumber \\
& & =  (I_{\phi }-\Delta_{\phi }) \big[  (I_{\phi }-\Delta_{\phi }  )^{(-1)}\circ {\rm d}J_{\phi }[\psi ][{\rm Tr}\, u_l]+ {\rm d}(  I_{\phi }-\Delta_{\phi }  )^{(-1)}[\psi ]\circ J_{\phi }[{\rm Tr}\, u_l]  \big] [u_m]\nonumber \\
& & =    {\rm d}J_{\phi }[\psi ][{\rm Tr}\, u_l][u_m]-  \lambda^{-1} {\rm d}(I_{\phi }-\Delta_{\phi }  )[\psi ][ u_l][u_m]  \, .
\end{eqnarray}

By standard calculus, by means of a change of variables and formula (\ref{der4}) it follows that

\begin{eqnarray}\label{dedelta1}
{\rm d}\Delta_{\phi }  [\psi ][u_l][u_m]= \int_{\phi (\Omega )}\nabla v_l(\nabla \zeta +\nabla \zeta^T  )\nabla v_m^Tdy -\int_{\phi (\Omega )} \nabla v_l\nabla v_m^T{\rm div }\zeta dy\, .
\end{eqnarray}

Since $\Delta v_l =v_l $ on $\phi (\Omega )$ and $\frac{\partial v_l}{\partial \nu }=\lambda v_l  $ on $\partial \phi (\Omega )$ for all $l\in F$, by the Divergence Theorem it follows that

\begin{eqnarray}\label{dedelta2}
\lefteqn{ \int_{\phi (\Omega )}\nabla v_l(\nabla \zeta +\nabla \zeta^T  )\nabla v_m^Tdy }\nonumber \\
& & =\sum_{r,s=1}^N\int_{\phi (\Omega )}\frac{\partial \zeta_r }{\partial y_s}\frac{\partial v_l }{\partial y_r}\frac{\partial v_m }{\partial y_s}+\frac{\partial \zeta_s }{\partial y_r}\frac{\partial v_l }{\partial y_r}\frac{\partial v_m }{\partial y_s}dy
  =\sum_{r,s=1}^N\int_{\partial \phi (\Omega )} \zeta_r\nu_s \frac{\partial v_l }{\partial y_r}\frac{\partial v_m }{\partial y_s}\nonumber \\
& & + \zeta_s\nu_r \frac{\partial v_l }{\partial y_r}\frac{\partial v_m }{\partial y_s}d\sigma  -
\sum_{r,s=1}^N \int_{\phi (\Omega )}\zeta _r\frac{\partial ^2v_l}{\partial y_r\partial y_s}\frac{\partial v_m}{\partial y_s}
+\zeta_r \frac{\partial v_l}{\partial y_r}\frac{\partial^2 v_m}{\partial y_s^2}
+\zeta_s \frac{\partial^2 v_l}{\partial y_r^2}\frac{\partial v_m}{\partial y_s}\nonumber \\
& & + \zeta _s\frac{\partial v_l}{\partial y_r}\frac{\partial ^2v_m}{\partial y_r\partial y_s}dy=\int_{\partial \phi (\Omega )}\biggl(\frac{\partial v_m}{\partial \nu }\nabla v_l +\frac{\partial v_l}{\partial \nu }\nabla v_m - (\nabla v_l\nabla v_m^T)\nu
 \biggr)\zeta^Td\sigma \nonumber \\
& & + \int_{\phi (\Omega )}{\rm div }\zeta \nabla v_l\nabla v_m^T-(\Delta v_m\nabla v_l+\Delta v_l\nabla v_m   )\zeta^Tdy =
\lambda \int_{\partial \phi (\Omega )}\left( v_m\nabla v_l \right. \nonumber \\
& & \left. + v_l\nabla v_m\right)\zeta^Td\sigma
+\int_{\phi (\Omega )}{\rm div }\zeta \nabla v_l\nabla v_m^T
-\int_{\phi (\Omega )}( v_m\nabla v_l+ v_l\nabla v_m   )\zeta^Tdy \nonumber \\
& & -\int_{\partial \phi (\Omega )}\nabla v_l\nabla v_m^T(\nu \zeta^T) d\sigma =\lambda \int_{\partial \phi (\Omega )} \nabla (v_lv_m)    \zeta^Td\sigma +\int_{\phi (\Omega )} ( \nabla v_l\nabla v_m^T \nonumber \\
& & + v_lv_m )\, {\rm div }\zeta dy -\int_{\partial \phi (\Omega )}(\nabla v_l\nabla v_m^T  +v_lv_m) (\nu \zeta^T) d\sigma .
\end{eqnarray}

Moreover, by (\ref{der4}) and a change of variable it follows that

\begin{equation}\label{dedelta3}
{\rm d}I_{\phi}[\psi ][v_l][v_m]=\int_{\phi (\Omega )}v_lv_m {\rm div }\zeta dy .
\end{equation}

By (\ref{dedelta1})-(\ref{dedelta3}) we obtain

\begin{eqnarray}\label{dedelta4}
{\rm d}(I_{\phi }-\Delta_{\phi })[\psi ][v_l][v_m]=
\int_{\partial \phi (\Omega )}[(\nabla v_l\nabla v_m^T  +v_lv_m) \nu -\lambda  \nabla (v_lv_m) ]\zeta^Td\sigma
\end{eqnarray}

which combined with (\ref{dedelta0}) and (\ref{techlem1}) yields

\begin{eqnarray}\label{dedelta5}
\lefteqn{<{\rm d}T_{\phi}[\psi][u_l], u_m>_{\phi}}\nonumber \\
& &=\lambda^{-1}\int_{\partial \phi (\Omega )}\left[\lambda\left( \frac{\partial (v_lv_m)}{\partial \nu} +Hv_lv_m\right)-(\nabla v_l\nabla v_m^T+v_lv_m )\right]\nu\zeta^Td\sigma \, .
\end{eqnarray}

By equalities
(\ref{sym4}), (\ref{sym5}) and (\ref{dedelta5})
it
follows that
\begin{eqnarray*}
\lefteqn{{\rm d}\Lambda_{F,h}[\phi][\psi]
 }
\\ \nonumber
& &
=\left\{
{\rm d} \Gamma _{F,|F|-h}[\phi][\psi]
\Gamma_{F,|F|} [\phi ]   \right.
-\left. \Gamma_{F,|F|-h} [\phi ]
{\rm d}\Gamma_{F,|F|} [\phi][\psi]
\right\}\lambda_{F}^{2|F|}[\phi ]
\\ \nonumber
& &
=\left[
{
|F|-1
\choose
|F|-h-1
}
\lambda_{F}^{h+1-2|F|}[\phi ]
-
\left(
\begin{array}{c}
|F|
\\
h
\end{array}
\right)
{
|F|-1
\choose
|F|-1
}
\lambda_{F}^{h+1-2|F|}[\phi ]
\right]
\\ \nonumber
& &
\cdot\lambda_{F}^{2|F|}[\phi ]
\sum_{l\in F}
<{\rm d}T_{\phi}[\psi][u_l], u_l>_{\phi}
=
-\lambda_{F}^{h+1}[\phi]
{
|F|-1
\choose
h-1
}
\sum_{l\in F}
<{\rm d}T_{\phi}[\psi][u_l], u_l>_{\phi}\\
& & =
\lambda_{F}^{h}[\phi]
{
|F|-1
\choose
h-1
}
\sum_{l\in F}
\int_{\partial \phi (\Omega )}\left[
|\nabla v_l|^2+v_l^2
-\lambda_F[\phi ]\left( \frac{\partial (v_l^2)}{\partial \nu} +Hv_l^2\right)
\right](\psi\circ \phi^{(-1)})\nu^T  d\sigma .
\end{eqnarray*}
Formula (\ref{sym2}) then  follows by observing that
$$|\nabla v_l|^2+v_l^2
-\lambda_F[\phi ] \left(\frac{\partial (v_l^2)}{\partial \nu} +Hv_l^2\right)=|\nabla_T v_l|^2+(1-\lambda_F[\phi]H-\lambda_F^2[\phi ])v_l^2,$$
for all $l\in F$.\hfill $\Box$\\

\section{Overdetermined problems}

In this section we consider isovolumetric and isoperimetric domain perturbations. Namely, given a bounded domain $\Omega$ in ${\mathbb{R}}^N$ of class $C^2$,
we consider domain transformations $\phi \in {\mathcal{A}}_{\Omega }$ satisfying either the volume constraint
\begin{equation}
\label{volume}
{\rm Vol }\, \phi (\Omega )={\rm const},
\end{equation}
or the perimeter constraint
\begin{equation}\label{perimeter}
{\rm  Per}\,  \phi (\Omega )={\rm const},
\end{equation}
where ${\rm Vol }\, \phi (\Omega )$ denotes the N-dimensional Lebesgue measure of $\phi (\Omega )$ and  ${\rm Per}\,  \phi (\Omega )$ the (N-1)-dimensional surface measure
of $\partial \phi (\Omega )$. It is here natural to introduce the functionals $V$ and $P$ from ${\mathcal {A}}_{\Omega }$ to ${\mathbb{R}}$ defined by
\begin{equation}
V[\phi ]={\rm Vol }\, \phi (\Omega ),\ \ \ {\rm and} \ \ \ P[\phi ]={\rm  Per}\,  \phi (\Omega ),
\end{equation}
for all $\phi \in {\mathcal{A}}_{\Omega }$. Note that by changing variables in integrals
\begin{equation}\label{volper}
V[\phi ]=\int_{\Omega }|{\rm det }\nabla \phi | dx,\ \ \ {\rm and } \ \ \ P[\phi ]=\int_{\partial \Omega } \left|\nu (\nabla \phi )^{-1}   \right| \left|\mathrm{det}\nabla \phi\right|d\sigma ,
\end{equation}
for all $\phi \in {\mathcal{A}}_{\Omega }$. We recall the following

\begin{defn} Let $\Omega$ be a bounded domain in ${\mathbb{R}}^N$ of class $C^2$.
Let ${\mathcal {F}}$ be a real-valued differentiable map defined
on an open subset of  ${\mathcal {A}}_{\Omega }$.
\begin{itemize}
\item[(i)] We say that $\phi\in {\mathcal{A}}_{\Omega }$
is a critical point for ${\mathcal{F}}$ with volume constraint
(\ref{volume}) if
\begin{equation}
\label{ker}
{\rm ker }\, {\rm d}V(\phi )\subseteq  {\rm ker }\, {\rm d}{\mathcal{F}}(\phi ) .
\end{equation}
\item[(ii)] We say that $\phi\in {\mathcal{A}}_{\Omega }$
is a critical point for ${\mathcal{F}}$ with perimeter  constraint
(\ref{perimeter}) if
\begin{equation}
\label{ker2}
{\rm ker }\, {\rm d}P(\phi )\subseteq  {\rm ker }\, {\rm d}{\mathcal{F}}(\phi ) .
\end{equation}
\end{itemize}
\end{defn}

As is well-known definitions (i) and (ii) are  related to the problem of finding  local extremal points for
$$
\min_{V[\phi ]={\rm const}}{\mathcal{F}}[\phi ]\ \ \ {\rm or}\ \ \ \max_{   V[\phi ]={\rm const}} {\mathcal{F}}[\phi],
$$
and
$$
\min_{P[\phi ]={\rm const}}{\mathcal{F}}[\phi ]\ \ \ {\rm or}\ \ \ \max_{   P[\phi ]={\rm const}} {\mathcal{F}}[\phi],
$$
respectively. Indeed if $\phi$ is a local minimum or maximum of a function ${\mathcal{F}}$ under condition (\ref{volume}) (resp. (\ref{perimeter}))   then inclusion (\ref{ker}) (resp. (\ref{ker2}))  holds (see, {\it e.g.}, Deimling~\cite[\S~26]{dei}).

Statements (i) and (ii) in the following theorem provide a characterization of all critical points of functions $\Lambda_{F,h}$ under volume or perimeter constraint. We note that the proof of statement (iii) in the following theorem makes use of the  same argument of \cite[Prop.~2.30]{lalacri} which is included here for the convenience of the reader.

\begin{thm}Let $\Omega $ be a bounded domain in ${\mathbb{R}}^N$ of class $C^2$ and $F$ be a finite subset of ${\mathbb{N}}$.
Assume that $\phi\in \Theta_{\Omega } [F]$ and that the eigenvalues $\lambda_j[\phi ]$ have the common value $\lambda_F[\phi ]$
for all $j\in F$. Let $\{ v_l\}_{l\in F}$ be an orthornormal basis in $W^{1,2}(\phi (\Omega ))$ of the eigenspace corresponding to $\lambda_F[\phi ]$.
\begin{itemize}
\item[(i)] The transformation $\phi$
is a critical point for any of the functions $\Lambda_{F,h}$, $h=1,\dots , |F|$,  with the volume constraint
(\ref{volume}) if and only if  there exists $C_1\in {\mathbb{R}}$ such that
\begin{equation}\label{ov1}
\sum_{l\in F}
 |\nabla_T v_l|^2+(1-\lambda_F[\phi]H-\lambda_F^2[\phi ])v_l^2
=C_1,\ \ {\rm on}\ \partial\phi (\Omega )\, .
\end{equation}
\item[(ii)] The transformation $\phi$
is a critical point for any of the functions $\Lambda_{F,h}$, $h=1,\dots , |F|$,  with the perimeter constraint
(\ref{perimeter}) if and only if  there exists $C_2\in {\mathbb{R}}$ such that
\begin{equation}\label{ov2}
\sum_{l\in F}
 |\nabla_T v_l|^2+(1-\lambda_F[\phi]H-\lambda_F^2[\phi ])v_l^2
=C_2H,\ \ {\rm on}\ \partial\phi (\Omega )\, .
\end{equation}
\item[(iii)] If $\phi (\Omega )$ is a ball then conditions (\ref{ov1}) and (\ref{ov2}) are satisfied.
\end{itemize}
\label{ov}
\end{thm}

{\bf Proof.} First we prove (i) and (ii). As for the case of finite dimensional vector spaces, the fact that $\phi  $ is a critical point for $\Lambda_{F,h}$
under volume or perimeter constraint respectively, is equivalent to the existence of the corresponding Lagrange multipliers which means that there exist
$c_1, c_2\in {\mathbb{R}}$ such that  ${\rm d}\Lambda_{F,h} (\phi  )=c_1 {\rm d}V(\phi  )$ and ${\rm d}\Lambda_{F,h} (\phi  )=c_2 {\rm d}P(\phi  )$  respectively (see, {\it e.g.}, Deimling~\cite[Thm.~26.1]{dei}).
By (\ref{techlem1}),  (\ref{der4}), (\ref{volper}), by changing variables in integrals and the Divergence Theorem it follows
that
\begin{equation}\label{dv}
{\rm d}V(\phi )[\psi ]= \int_{\partial \phi (\Omega )}(\psi\circ\phi ^{(-1)})\nu^Td\sigma ,\ \ {\rm and}\ \ {\rm d}P(\phi )[\psi]= \int_{\partial \phi (\Omega )}H(\psi\circ\phi ^{(-1)})\nu^Td\sigma\, ,
\end{equation}
for all $\psi \in C^2(\bar \Omega\, ; {\mathbb{R}}^N)$. Thus by combining  (\ref{sym2}), (\ref{dv}) and by the arbitrary choice of $\psi$ we deduce the validity of (\ref{ov1}) and (\ref{ov2}) for suitable constants $C_1, C_2$.

In order to prove (iii) we  proceed as in \cite{lalacri}. Without any loss of generality, we assume that $\phi (\Omega )$ is a ball centered at zero.  By rotation
invariance of the Laplace operator,  $\{v_{l}\circ A\}_{l\in F}$   is an orthonormal basis of the eigenspace  corresponding to $\lambda_F[\phi ]  $ for all $A\in O_{N}({\mathbb{R}})$. Here $O_{N}({\mathbb{R}})$ denotes the group of orthogonal linear transformations in ${\mathbb{R}}^N$.
 Since both $\{v_l\}_{l\in F}$
and $\{ v_{l}\circ A\}_{l\in F}$ are orthonormal bases of the same space, it follows that
\begin{equation}
\label{balld1}
\sum_{l=1}^{|F|} {v}_{l}^{2}\circ
A=\sum_{l=1}^{|F|} {v}_{l}^{2}\, ,
\end{equation}
 for all $A\in O_{n}({\mathbb{R}})$. Thus $\sum_{l=1}^{|F|}{v}_{l}^{2}$ is radial. This, combined with the fact that in the case of a ball $H$ is constant, implies the validity of  (\ref{ov1}) and (\ref{ov2}) for suitable constants $C_1$, $C_2$.\hfill $\Box$\\

Theorem~\ref{ov} suggests to consider the following overdetermined systems where $\lambda $ is an eigenvalue of multiplicity $m$ and $u_1, \dots , u_m$ is an orthonormal basis in $W^{1,2}(\Omega )$ of the corresponding eigenspace

\begin{equation}\label{pr1}
{\rm (P1)}\left\{\begin{array}{ll}
\Delta u_l=u_l,\ \  & {\rm in }\ \Omega \, ,\vspace{1mm}\\
\frac{\partial u_l}{\partial \nu}=\lambda u_l, \ \  & {\rm in }\ \partial \Omega  \, ,\vspace{1mm} \\
l=1,\dots , m, & \\
\sum_{l=1}^m
|\nabla_T v_l|^2+(1-\lambda H-\lambda ^2)v_l^2       ={\rm const},\ \ &  {\rm on}\ \partial\Omega \, ,
\end{array}
\right.
\end{equation}

\begin{equation}\label{pr2}
{\rm (P2)}\left\{\begin{array}{ll}
\Delta u_l=u_l,\ \ & {\rm in }\ \Omega \, , \vspace{1mm}\\
\frac{\partial u_l}{\partial \nu}=\lambda u_l, \ \ & {\rm in }\ \partial \Omega \, , \vspace{1mm} \\
l=1,\dots , m, & \\
\sum_{l=1}^m
 |\nabla_T v_l|^2+(1-\lambda H-\lambda ^2)v_l^2        ={\rm const}   H,\ \ &  {\rm on}\ \partial\Omega \, .
\end{array}
\right.
\end{equation}

It would be interesting to characterize those bounded domains $\Omega $ in ${\mathbb{R}}^N$ of class $C^2$  such that systems {\rm (P1)} or {\rm (P2)} are satisfied and to know whether the existence of solutions to one of these systems implies that $\Omega $ is a ball.\\

\vspace{0.5cm}

\noindent
{\small Pier Domenico Lamberti\\
Dipartimento di Matematica Pura ed Applicata\\
Universit\`{a} degli Studi di Padova\\
Via Trieste 63\\
35121 Padova\\
Italy\\
e-mail: lamberti@math.unipd.it}

\end{document}